\documentclass[12pt]{article}
 
\makeatletter
\@addtoreset{equation}{section}

\makeatother

\usepackage{amsfonts}

\def\ea{{\it et al.}\ }

\def\R{{\mathbb R}}
\def\Z{{\mathbb Z}}
\def\N{{\mathbb N}}

\def\PN{{\cal N}}
\def\PM{{\cal M}}

\def\cD{{\cal D}}

\def\sqr{\vcenter{
         \hrule height.1mm
         \hbox{\vrule width.1mm height2.2mm\kern2.18mm\vrule width.1mm}
         \hrule height.1mm}}                  
\def\square{\ifmmode\sqr\else{$\sqr$}\fi}
\def\qed{\hfill \square}

\def\be{\begin{equation}}
\def\ee{\end{equation}}

\def\eps{\varepsilon}
\let\phi=\varphi
\def\DD{\displaystyle}
\def\Ind#1{{\bf 1}_{\{#1\}}}
\def\P{{\bf P}}
\def\E{{\bf E}}
\def\Ee#1{{\bf E}^e_{#1}}
\def\Ev{{\bf E}^v}
\def\Eh#1#2{{\bf E}^h_{#1,#2}}
\def\Pe#1{{\bf P}^e_{#1}}
\def\Pv{{\bf P}^v}
\def\Ph#1#2{{\bf P}^h_{#1,#2}}
\def\rh#1{\rho\big(#1\big)}
\def\zz{\{0,1\}^{\Z}}
\def\epss{\eps^{\prime\prime}}

\newtheorem{theo}{Theorem}[section]
\newtheorem{lmm}{Lemma}[section]
\newtheorem{df}{Definition}[section]
\newtheorem{rem}{Remark}[section]
\newtheorem{conj}{Conjecture}[section]

\title{A Mixture of the Exclusion Process and the Voter Model}
\author{Vladimir~Belitsky$^1$ \and Pablo~A.~Ferrari$^1$ \and
Mikhail~V.~Menshikov$^{*,1,2}$ \and Serguei~Yu.~Popov$^{1,3}$}

\date{}

\begin{document}

\maketitle

{\small
\noindent
$^1$ and the mailing address: 
Instituto de Matem\'atica e Estat\'\i{}stica, Universidade de S\~ao Paulo, 
Caixa Postal 66281, 05315--970 - S\~{a}o Paulo, SP -
    BRAZIL

\smallskip
\noindent
$^2$ Faculty of Mathematics and Mechanics, Moscow
State University

\smallskip
\noindent
$^3$ Institute for Problems of Information Transmission,
Russian Academy of Sciences

\smallskip
\noindent
E-mails: vbel1@ime.usp.br, pablo@ime.usp.br, menchik@ime.usp.br,
popov@ime.usp.br}

\begin{abstract}
  We consider a one-dimensional nearest-neighbor interacting particle
  system, which is a mixture of the simple exclusion process
  and the voter model. The state space is taken to be the countable set
  of the configurations that have a finite number of particles to the right
  of the origin and a finite number of empty sites to the left of it.
  We obtain criteria for the ergodicity and some other properties of
  this system using the method of Lyapunov functions.

\end{abstract}



\noindent{\bf Running title:} A mixture of Exclusion Process and Voter Model

\section{Introduction}
\label{intro}
In this paper we consider a process that is a mixture of two nearest-neighbor 
one-dimensional interacting particle systems: the simple exclusion process
and the voter model. Let us first define these two processes
\begin{df}
\label{df1}
For $\eta\in\zz$ denote
 $$ \eta_{x,y}(z)  = 
    \left\{
  \begin{array}{ll}
   \eta(y), & \mbox{if } z=x,\\
   \eta(x), & \mbox{if } z=y,\\
   \eta(z), & \mbox{if } z\neq x,y,
  \end{array}
\right.$$ 
and
 $$ \eta_{x}(z)  = 
    \left\{
  \begin{array}{ll}
   1-\eta(z), & \mbox{if } z=x,\\
   \eta(z), & \mbox{if } z\neq x.
  \end{array}
\right.$$
A Markov process $\eta_t \in \zz$, $t\in [0, +\infty)$ is called
\begin{itemize}
\item Simple exclusion process with parameter $0\leq
  p \leq 1$, if its generator $\Omega_p^e$ has the following
  form:
\begin{eqnarray*}
 \Omega_p^e f(\eta) &=& 
\sum_{x,y} p(x,y)\,\eta(x)(1-\eta(y)) 
       [f(\eta_{x,y})-f(\eta)],
\end{eqnarray*}
where 
$$
  p(x,y) = \left\{\begin{array}{ll}
p, &\hbox{if }y=x-1,\\
1-p,&\hbox{if }y=x+1,\\
0,&\hbox{otherwise;}\end {array}\right.
$$
\item voter model, if its generator $\Omega^v$ is defined in the
  following way:
 $$ \Omega^v f(\eta) = \sum_x c(x,\eta)[f(\eta_x)-f(\eta)], $$
where 
\be
\label{vovochka}
 c(x,\eta)  = 
    \left\{
  \begin{array}{ll}
  \DD \frac12 (\eta(x-1)+\eta(x+1)), & \mbox{if } \eta(x)=0,\\  
\\ \DD \frac12 (2-\eta(x-1)-\eta(x+1)), & \mbox{if } \eta(x)=1.
  \end{array}
\right.
\ee 
\end{itemize}

\end{df}

The construction of those processes from their generators may be found in
the book of Liggett (1985); see the first chapter and the beginning of
the chapters corresponding to those processes. Harris graphical
construction (see Durrett (1988, 1995) for instance) is an alternative
approach to define these processes. It will be briefly reviewed and
used in Section~\ref{discrete}.

Let us call $\eta\in \zz$ {\it a configuration of particles} and let
us interpret $\eta(x)=1$ as the presence of a particle at the site
$x\in\Z$ in the configuration~$\eta$ and $\eta(x)=0$ as the absence of
it. The dynamics of both processes may be interpreted in terms
of particles that hop on~$\Z$ (the case of the exclusion process) or
appear and disappear at the sites of~$\Z$ (the case of the voter
model).

In the exclusion process, there may be at most one particle at each
site of~$\Z$. If there is a particle at site $x$ and no particle at
site $x+1$ (respectively at site $x-1$), then the particle at $x$
jumps with rate $(1-p)$ (respectively $p$) to site $x+1$ (respectively
$x-1$). This is a conservative dynamics, in the sense that neither
particles are created nor disappear. Liggett (1976) described the set
of invariant measures for this process. If $p=1/2$, the invariant
measures are convex combinations of the translation invariant product
measures parameterized with the density of particles. If $p> 1/2$, the
set of invariant measures contains also measures with support in the
countable state space
\begin{eqnarray}
  \label{eq:p1}
  \cD\, :=\, \hbox{the set of configurations with a
finite number of empty sites to the }\nonumber\\
\hbox{left of the origin and a finite
number of particles to the right of it.}\nonumber 
\end{eqnarray}
These measures are called \emph{blocking measures} because, due to the
exclusion rule and the accumulation of particles to the left of the
origin, the flux of particles is null. Of course there are also
blocking measures for $p<1/2$; they are obtained by the reflection
($\Z\rightarrow -\Z$) of those mentioned above. 
When an asymmetric exclusion process ($p\neq1/2$) is considered
from a random position determined by a so-called second-class particle,
a new set of invariant measures arises. They are called \emph{shock
  measures}, they have support on configurations with different
asymptotic densities to the left and right of the origin.
The respective results have the origin in the works of
 Ferrari \ea
(1991) and
Ferrari (1992). See the review paper of Ferrari (1994) and the book of
Liggett (1999) for an account of properties of these measures and the
asymptotic behavior of the second class particle.  Derrida \ea
(1998) propose a nice alternative descriptions of shock measures for
this process.

In the voter model, there may be at most
one particle per site, however its dynamics is nonconservative: 
 a new particle is born at an empty site $x$
at a rate proportional to the number of nearest neighbors of $x$
occupied by particles; and a particle that is present at a site $x$
disappears at a rate proportional to the number of the empty neighbors
of $x$.  Since only one site changes its value at any given time, this
model is a particular case of the so called \emph{spin-flip} models.
There are only two invariant measures for the one-dimensional voter
model defined above: one has the support on the configuration ``all
zeros'' and the other one has the support on the configuration ``all
ones''. The basic tool to prove those results is \emph{duality}, a
technique that allows to express properties of the voter model as
properties of a dual process,
 a process obtained when one ``looks backwards in
time''. There are two dual processes for the voter model: coalescing
random walks and annihilating random walks. See Liggett (1985,
Chapter~V), Durrett (1995) for accounts on these and many other
properties of the voter model.

If the voter model starts from the Heaviside configuration $\eta^0$,
defined by $\eta^0(x) = \Ind{x\le 0}$, then at any future time it is
a random translation of~$\eta^0$.  Indeed, the position of the
rightmost particle $X_t=\max \{x: \eta_t(x)=1\}$ performs a nearest
neighbor symmetric random walk and $\theta_{X_t}\eta_t = \eta_0$,
where $\theta_x$ is translation by $x$. This example motivates the
introduction of an equivalence relation: we say that two
configurations $\eta$ and $\eta'$ are equivalent and write
$\eta\sim\eta'$ if one of them is a translation of the other: there
exists a $y\in \Z$ such that $\eta(x)=\eta'(x+y)$ for all $x\in\Z$.
Let $\tilde \cD:=\cD/\sim$ denote the set of equivalence classes
induced by $\sim$.  Let then $\cD_0$ denote the set of the
configurations in the equivalence class of $\eta^0$. In the voter
model, $\eta_0\in\cD_0$ implies $\eta_t\in\cD_0$ for all $t$.  Hence,
denoting $\tilde \eta_t$ the equivalence class of $\eta_t$, we have
that $\tilde\eta_0=\tilde\eta^0$ implies $\tilde\eta_t\equiv
\tilde\eta^0$ (nothing moves). The process $\tilde\eta_t\in\tilde\cD$
just defined, is isomorphic to $\theta_{X_t}\eta_t$, the voter model
as seen from its rightmost particle.

Cox and Durrett (1995) studied one dimensional voter models on $\tilde
\cD$ with rate function $\sum_y q(|x-y|)
|\eta(x)-\eta(y)|$ for some probability function $q(x)$.  They show
that if $\sum_x|x|^3q(x) < \infty$, then the process as seen from the
rightmost particle $\tilde\eta_t\in\tilde\cD$ is positive recurrent
and hence admits a unique invariant (shock) measure. Calling $Y_t$ the
leftmost hole, this implies that under the invariant measure the size
of the \emph{hybrid zone} ---the region of coexistence of zeros and
ones--- $X_t-Y_t$ is bigger than $-1$ and finite with probability one;
and of course its distribution is independent of $t$. They also prove
that the expected value of $X_t-Y_t$ under the invariant measure is
infinite and that $X_t/\sqrt t$ converges as $t\to\infty$ to a
centered normal distribution with finite variance.  The approach is
based on a fine analysis of the (dual) process coalescing random
walks.  It is also shown there that there are no ``stable'' hybrid
zones in dimension $d=2$: if one starts with ones in the negative $x$
semiplane and zeros in the positive semiplane and paints 1s white and
0s black, then the normal distribution with variance $\sim t$ predicts
the shade of grey we see at time $t$ in the horizontal direction.

Ferrari (1996) shows the existence of an invariant shock measure for
the biased voter model as seen from the rightmost particle. In this
model the rate function is given by $c_2(x,\eta) = (a \eta(x) + b(1-\eta(x)))
c(x,\eta)$, with $c(x,\eta)$ as defined in~(\ref{vovochka}). The
proof in this case is more direct because it is based on
straightforward dominations by supermartigales.

The goal of this paper is the study of the existence of shock measures
in a {\it mixture\/} of the exclusion process and the voter model.

\begin{df}
\label{df2}
Let $\beta\in[0,1]$. A Markov process $\eta_t \in \zz$, $t\in [0,
+\infty)$ is called \emph{hybrid process} with mixing
parameter~$\beta$ and exclusion parameter~$p$, if its generator is
\be
  \label{h1}
  \Omega^h_{\beta,p}:= (1-\beta)\Omega^e_p + \beta \Omega^v.
\ee

\end{df}

The hybrid process $\tilde\eta_t$ (the class of equivalence of
$\eta_t$ with initial configuration in $\cD$) is a Markov process on
$\tilde\cD$. This process is a particular case of a model of random
grammars, considered by Malyshev (1998). Models consisting of a
mixture of a spin-flip dynamics and a symmetric exclusion dynamics are
usually called in the literature ``diffusion-reaction processes''.
When $\beta \gg 0$, an appropriate space-time rescaling with $\beta$
produces hydrodynamic limits rising the reaction-diffusion equation
${\partial u \over \partial t} = {\partial ^2 u\over \partial^2 x} +
f(u)$, whereas the function $f$ is related to the spin-flip dynamics
and $u=u(x,t)\in[0,1]$, $x,t\in \R_+$ corresponds to the macroscopic
density of particles (De Masi \ea
(1986)). In some
cases these equations accept traveling-wave solutions ---solutions of
the type $u(x,t) = u_0(x-vt)$ for some speed $v$ with
$\lim_{x\to\infty} u_0(x) = 0$, $\lim_{x\to-\infty} u_0(x) = 1$. This
motivates the question about the existence of a microscopic
counterpart of the macroscopic traveling wave solutions.  A particular
case of reaction process is the growth model, a process with rate
function $c(x,\eta)(1-\eta(x))$, where $c(x,\eta)$ has been 
defined in~(\ref{vovochka}): 0 flips to 1 at rate proportional
to the number of ones in the neighborhood, but 1 never flips to 0.
Bramson \ea
(1986)
showed the existence of an invariant (shock) measure for the process
$\tilde\eta_t$, where $\eta_t$ is any nontrivial mixture of the
exclusion process and the growth model. Cammarota and Ferrari (1991)
proved the Normal asymptotic behavior of $(X_t-\E X_t)/\sqrt t$ for
this mixture.  Machado (1998) studied this process in a strip and in
$\Z^d$.

Let $\tilde\tau_c(\tilde\eta)$ be the first time the process $\tilde\eta_t$
starting with the configuration $\tilde\eta\in\tilde\cD$ hits
$\tilde\eta^0$, the Heaviside configuration defined above. The
subscript $c$ refers to continuous time (as a counterpart of a
discrete-time process to be introduced below).  Let us recall some classical
definitions. We say that the process $\tilde\eta_t$ is {\it
  transient}, if $\P(\tilde\tau_c(\tilde\eta) <\infty)<1$ and {\it
  recurrent}, if $\P(\tilde\tau_c(\tilde\eta) <\infty)=1$. In the last case
we say that the process is {\it positive recurrent} if
$\E(\tilde\tau_c(\tilde\eta))<\infty$ and {\it null recurrent} if this
expectation is infinity. An irreducible countable Markov chain is
\emph{ergodic} if it has a unique invariant measure. Since,
except for the pure voter model,
$\tilde\eta_t$ is irreducible, positive recurrence is equivalent to
ergodicity in our context. The following theorem contains our results.

\begin{theo}\label {pp}
  Let $\eta_t$ be a process in $\cD$ with generator
  $\Omega^h_{\beta,p}$. Let $\tilde\eta_t$ be the corresponding
  process in the space of classes of equivalence $\tilde\cD$.
\begin{enumerate}
\item Exclusion process. Assume $\beta = 0$. Then the process
  $\tilde\eta_t$ is ergodic for $p>1/2$ and transient for $p\le 1/2$.
  
\item Hybrid process. Assume $0<\beta<1$. Then 
\subitem {\rm i)} There exists
  $\beta_c<1$ such that for any $\beta>\beta_c$ and any $p\in (0,1)$
  the process $\tilde\eta_t$ is ergodic.  
\subitem {\rm ii)} For any $p\ge 1/2$ and any $\beta$,
  the process is ergodic.
  
\item Voter model. Assume $\beta=1$. Then the process $\tilde\eta_t$
  is positive recurrent.  Moreover, for any initial configuration
  $\tilde\eta\in\tilde\cD$ and any $\eps>0$, \be
\label{p2}
\E (\tilde\tau_c(\tilde\eta))^{3/2-\eps} < \infty\,;\;\;\;\E
(\tilde\tau_c(\tilde\eta))^{3/2+\eps} = \infty.  \ee
\end{enumerate}
\end{theo}

The fact that the exclusion process $\tilde\eta_t$ in $\tilde\cD$ is
ergodic for $p>1/2$ follows immediately from well known results of
Liggett (1976, 1985) who described the invariant measures for
$\eta_t$ in the irreducible classes of $\cD$. Since the system is
conservative, ergodicity of $\eta_t$ on any irreducible class of $\cD$
is equivalent to ergodicity of $\tilde\eta_t$ on $\tilde\cD$. Our
alternative approach does not use the knowledge of the invariant
measure. When $p\le 1/2$, the results of Liggett imply only that the
process $\eta_t$ is not positive recurrent; our result says that it is
transient. For $p<1/2$, the transience holds immediately from laws of large
numbers for the leftmost hole and the rightmost particle. For $p=1/2$,
the transience is a more delicate matter.

The bounds in (\ref{p2}) show the velocity of the convergence of the
voter model to the invariant measure, which is the singleton supported
by $\cD_0$. It may be the case that these bounds could be obtained
from the duality of the voter model to the coalescing random walks,
however, we have not investigated this approach.

Our main results are the conditions for ergodicity for the hybrid
model described in point (2) of the theorem. It says that if either
the proportion of voter in the hybrid process is large enough or the
exclusion process has no drift to the right, then the hybrid process
is ergodic. Item (1) says that exclusion is transient for $p\le 1/2$,
while item (3) says that voter is always 
positive recurrent. The first part of
item (2) says that voter ``wins'' if the proportion of voter is
sufficiently large, uniformly on the exclusion asymmetry; the relevant
point in the second part of item (2) says that for the symmetric
exclusion, any proportion of voter guarantees ergodicity.

We are not totally satisfied with this result because
sufficient conditions for transience are missing. One would like to
show that if the asymmetry of the exclusion process has a tendency to
``escape'' from $\cD$ then an addition of a small proportion of the
voter model will not be able to prevent it from escaping. But for
now, it is still very unclear to us, if the process could be
transient in this case.  We state
now a conjecture for the nonergodicity of the hybrid process. A
heuristic argument supporting the conjecture is presented in
Section~\ref{hybrid}. 

\begin{conj}
\label{pcnj} For any $p<1/2$ there exists a $\beta_0(p)>0$ such that
for any $\beta<\beta_0(p)$, the hybrid process $\tilde\eta_t$ with
parameters $\beta$ and $p$ is not ergodic.
\end{conj}

The parameter space $\{(p,\beta)\,:\, p, \beta\in[0,1]\}$ is
partitioned in three regions: ergodicity, transience and
null-recurrence. Presumably the region of transience satisfies the
property: if the hybrid process with parameters $(p_0,\beta_0)$ is
transient, then the one with parameters $(p_1,\beta_1)$ will be also
transient for $p_1\le p_0$ and $\beta_1\le\beta_0$. But we do not have
any monotonicity argument at hand to argue this. We know that the
transience region is nonempty because it contains the segment
$[0,1/2]\times\{0\}$, but we do not know how to prove that it contains
points in the interior of the parameter space.

How stable under changes of the dynamics are our results? Can we
extend Theorem~\ref{pp} to nonnearest-neighbors processes? When
$\beta=1$, only voter, the answer is given by Cox and Durrett (1995),
as described above. When $\beta=0$, only exclusion, it is known that
the process is not ergodic on $\cD$ if $p(x,y)$ is symmetric (all
invariant measures are translation invariant in this case), but it
is an open problem of Liggett (1985, Section~VIII.7, Problem~6)
 in the case when $p(x,y)$ is
asymmetric. The conjecture is that if $p(x,y) = q(y-x)$, for some $q$,
then the system would be ergodic under the condition $\sum_x x q(x)
<0$. In the final remarks we explain where our approach fails to
work when extended to the nonnearest-neighbors case.

Motivations coming from real life, description of shock measures in
other one-dimensional models and nice conjectures about the existence
of shock measures in other systems can be found in the introduction
of Cox and Durrett (1995).

Theorem~\ref{pp} is proven for the discrete-time version of
$\tilde\eta_t$ and then standard arguments are used to prove the
continuous counterpart. The discrete process is a Markov chain in
$\tilde\cD$. The basic tool is a set of theorems from Fayolle \ea
(1995), which give conditions for ergodicity, recurrence and
transience of denumerable Markov chains using so-called Lyapunov functions.
The application of these functions to the processes in interest
 produces sub or super martingales, which can be used
straightforwardly to show the desired properties. The problem is that
these functions are frequently hard to find. One of the contributions
of this paper is the exhibition of Lyapunov functions that work
for the exclusion process, the voter model and their mixture.
 
The paper is organized in the following manner. In Section~\ref{discrete} 
we introduce the discrete version of the process $\tilde\eta_t$.
In Section~\ref{pfmm} we state the
results of Fayolle \ea
(1995) we need. In Section~\ref{notation} we introduce the Lyapunov
functions of the process that will be relevant in the proofs. In
Sections~\ref{excl}, \ref{voter} and~\ref{hybrid} we state and prove
the results for the discrete-time versions of the exclusion process,
the voter model and the hybrid process respectively. 
In Section~\ref{continuous} we show how to pass from the discrete to the
continuous time and prove Theorem~\ref{pp}.

\section{Discrete and continuous-time processes}\label{discrete}

In this section we introduce discrete-time versions of the exclusion 
process, voter model, and their mixture that have been
 defined in the previous section, and establish their
relations with the continuous-time processes.

Let $\eta$ be a configuration from $\zz$.  We say that a discrepancy
of type $01$ ($10$) occurs in $\eta$ at the site $x$, if $\eta(x-1)=0,
\eta(x)=1$ (resp., $\eta(x)=1, \eta(x+1)=0$). The above defined
countable set $\cD$ is the set of those configurations of $\zz$ in
which there is only a finite number of discrepancies, and the number
of discrepancies of type $10$ minus the number of discrepancies of
type $01$ is equal to~$1$.  Then it is easy to see that
\begin{eqnarray*} 
  \cD & = & \{ \eta\in \zz : \mbox{there exist } i_0, j_0
 \mbox{ such that} \\
     && ~~~~ \eta(i) = 1 \mbox{ for } i\leq i_0
         \mbox{ and } \eta(j) = 0 \mbox{ for } j\geq j_0 \}  
\end{eqnarray*}
and that $\cD$ is countable.

The discrete time exclusion process with parameter $p$ (to be called
here EP($p$)) is a Markov process with the state space $\cD$ and the
following dynamics: for every $n\geq 0$, if $\eta$ is the state at
time $n$ then $\eta^\prime$, the state at time $n+1$, is obtained by
the following procedure (i)--(ii):
\begin{itemize}
\item[{\rm (i)}] we choose one of the discrepancies of $\eta$
with uniform distribution; say the discrepancy
 at the site $x$ has been chosen,
then  
\item[{\rm (ii)}] if the discrepancy is $01$ ($10$)
then we exchange $0$ and $1$ with the
 probability $0<p<1$ (resp., $0<q:=1-p<1$)
while nothing is changed with the resting probability $q$
 (resp., $1-q$).
\end{itemize}
The Exclusion Process just defined is a countable Markov chain on
$\cD$.

Let us define now the discrete time Voter Model (to be called VM)
and the discrete time hybrid process (to be called HP($\beta$, $p$),
where $\beta$ is the mixing parameter and $p$ is the
exclusion parameter).
For VM the step (i) is the same, and (ii) is 
substituted by the following:
\begin{itemize}
\item[{\rm (ii$^\prime$)}] the chosen discrepancy is substituted
by either $11$ or $00$
with probabilities $1/2$.
\end{itemize}
To construct HP($\beta$, $p$), we first execute (i), and then with
probability $1-\beta$ we execute (ii) (i.e. make a step of the
exclusion process), and with probability~$\beta$ we execute
(ii$^\prime$) (i.e. make a step of the voter model). We use the
notation $(\xi_n\,:\,n\in \N)$ for the HP($\beta$, $p$). $\xi_n$
denotes the configuration of the system at time $n$.
  
In~(\ref{h5}) below we shall present the relation between
the discrete-time hybrid process $(\xi_n\,:\,n\in \N)$ and the
continuous-time hybrid process $(\eta_t\,:\,t\geq 0)$ with the
mixing parameter~$\beta$ and the exclusion parameter~$p$.
To this end, we shall need the Harris graphical construction for
$(\eta_t\,:\,t\geq 0)$, which we now briefly recall.
It is a ``superposition'' of the graphical construction
for the voter model (see Durrett (1995)) with that for the exclusion
process (see Ferrari (1992)) with the respective weights~$\beta$
and~$(1-\beta)$.

Let $\{(\PN^{x,x+1}_t, t\geq 0)\}_{x\in\Z}$,
$\{(\PN^{x,x-1}_t, t\geq 0)\}_{x\in\Z}$, 
$\{(\PM^{x,x+1}_t, t\geq 0)\}_{x\in\Z}$,
$\{(\PM^{x,x-1}_t, t\geq 0)\}_{x\in\Z}$ be four independent families of
Poisson point processes 
 with the respective rates $(1-\beta) p$, $(1-\beta) q$,
$\beta /2$ and $\beta /2$. 
Given the initial configuration~$\eta_0$, the dynamics of
the process $\eta_t, t\geq 0$ is determined by those Poisson
processes in the following manner.
If there is a Poisson event at time~$t$ in $\PN^{x,x+1}$ 
(resp., $\PN^{x,x-1}$), which means $\PN^{x,x+1}_t-\PN^{x,x+1}_{t^-}=1$,
and if~$x$ has a particle while~$x+1$ is empty (resp., $x-1$ is empty)
in $\eta_{t^-}$, then the particle jumps from~$x$ to~$x+1$
(resp., $x-1$) at time~$t$. 
If there is a Poisson event at time~$t$ in $\PM^{x,x+1}$ 
(resp., $\PM^{x,x-1}$), then the site~$x+1$ (resp., $x-1$)
acquires the same state at time~$t$ as the state of~$x$ in $\eta_{t^-}$.

 Let $\tau_0=0$ and
for $n\ge 1$, set
\begin{eqnarray}
  \label{h2}
  \tau_n&=&\inf\Bigl\{t>\tau_{n-1}:\sum_{x,y:|x-y|=1}
  |\eta_{t^-}(x)-\eta_{t^-}(y)|\nonumber\\
  &&\qquad\qquad\qquad
  \times (\PN^{x,y}(\tau_{n-1},t]+\PM^{x,y}(\tau_{n-1},t])
  \;>\;0\Bigr\}    
\end{eqnarray}
where $\PN(s,t]$ denotes the number of the
 Poisson events in the time interval $(s,t]$
for the process $\PN$.  
We call $\tau_n$ the instants
of \emph{attempted jumps} of the process $\eta_t$. 
 It follows then from our
definitions that if $\eta_0=\xi_0$, then
\be
  \label{h5}
  (\xi_n\,:\,n\ge 0)\; = \;(\eta_{\tau_n}\,:\,n\ge 0)\qquad\hbox{in
  distribution.} 
\ee

\section{Criteria for recurrence and transience of Markov chains}\label{pfmm}

In this section we state the criteria for ergodicity, recurrence and
transience of countable Markov chains to be used in the sequel.
The next four theorems are Theorems 2.2.3, 2.2.1,  2.2.2, 2.2.7, respectively,
of Fayolle \ea
(1995).

\begin{theo}
\label{cr_erg}
Let $\xi_t$, $t=0,1,2,\ldots$ be an irreducible Markov chain with
the countable state space~$X$. Suppose that there exist a positive
function $f(x)$ and a finite set $A\subset X$ such that
\be
\label{eq_erg}
\E (f(\xi_{t+1})-f(\xi_t) \mid \xi_t = x ) \leq - \eps
\ee
for some $\eps>0$ and all $x\in X\setminus A$, and that
\be
\label{eq_erg1}
\E (f(\xi_{t+1}) \mid \xi_t = x ) < \infty
\ee
for $x\in A$. Then
the Markov chain is ergodic.
\end{theo}

\begin{theo}
\label{cr_rec}
Let $\xi_t$, $t=0,1,2,\ldots$ be an irreducible Markov chain with
the countable state space~$X$. Suppose that there exist a positive
function $f(x)$, $f(x)\to\infty$ as $x\to\infty$,
 and a finite set $A\subset X$ such that
\be
\label{eq_rec}
\E (f(\xi_{t+1})-f(\xi_t) \mid \xi_t = x ) \leq 0
\ee
for all $x\in X\setminus A$. Then
the Markov chain is recurrent.
\end{theo}

\begin{theo}
\label{cr_tr1}
Let $\xi_t$, $t=0,1,2,\ldots$ be an irreducible Markov chain with
the countable state space~$X$. Suppose that there exist a positive
function $f(x)$ 
 and a set $A\subset X$ such that~(\ref{eq_rec})
holds for all $x\in X\setminus A$ and
$$f(x_0) < \inf_{x\in A} f(x) $$
for some $x_0\notin A$. Then
the Markov chain is transient.
\end{theo}

\begin{theo}
\label{cr_tr2}
Let $\xi_t$, $t=0,1,2,\ldots$ be an irreducible Markov chain with
the countable state space~$X$. Suppose that there exist a positive
function $f(x)$ and a constant~$C$ such that
if $f(x)>C$, then
\be
\label{eq_tr2}
\E (f(\xi_{t+1})-f(\xi_t) \mid \xi_t = x ) \geq  \eps
\ee
for some $\eps>0$,
and suppose that for some $K>0$
\be
\label{eq_tr2'}
|f(\xi_{t+1})-f(\xi_t)| \leq K \quad \mbox{a.s.} 
\ee
 Then the Markov chain is transient.
\end{theo}

Besides the ergodicity, we are going to study the
existence of moments of the hitting time of the set~$\cD_0$. To do this, we
shall need the following result of Aspandiiarov \ea
(1996, Theorem~1)

\begin{theo}
\label{cr_mom}
 Let $A$ be some positive real number. Suppose that we
are given a $\{{\cal F}_n\}$-adapted stochastic process~$X_n$,
$n\geq 0$, taking values in an unbounded subset of $\R_+$.
Denote by $\tau_A $ the moment when the process~$X_n$
enters the set $(0,A)$.
Assume that there exist $\lambda>0$, $p_0\geq 1$ such that
for any $n$, $X_n^{2p_0}$ is integrable and
\be
\label{iasno}
 \E (X_{n+1}^{2p_0} - X_n^{2p_0} \mid {\cal F}_n) \leq
     \lambda X_n^{2p_0-2}
\ee
on $\{\tau_A>n\}$. Then there exists a positive constant
$C=C(\lambda,p_0)$ such that for all $x\geq 0$ whenever $X_0=x$
with probability~$1$
\be
\label{exists!!!}
\E\tau_A^{p_0} \leq C x^{2p_0}.
\ee
\end{theo}

\section{Functions of the process}
\label{notation}
For the sake of brevity we will substitute in the sequel
the expression ``block of zeros'' by ``$0$-block'' and
``block of ones'' by ``$1$-block''.

A class of equivalence $S\in \tilde\cD$ can be identified by a finite
set of positive numbers in the following form: 
\be
\label{m_in_i}
S=\ldots 111
  \overbrace{0000}^{n_1}
  \overbrace{11111}^{m_1}
  \overbrace{0000}^{n_2}
  \overbrace{11111}^{m_2}\ldots
  \overbrace{00000}^{n_N}
  \overbrace{1111}^{m_N}000\ldots ,
\ee
where $n_i = n_i(S)$ is the size of $i$-th $0$-block, $m_i = m_i(S)$
is the size of $i$-th $1$-block, $N = N(S)$ is the number of
$1$-blocks not including the leftmost infinite $1$-block.  In the
sequel the word ``configuration'' will usually mean ``class of
equivalence''.  So, for $S\in \tilde\cD$ we can simply write $S =
(n_1,m_1,\ldots,n_N, m_N) $. 

Denote
$r_0=0$, $r_i=\sum_{j=1}^i(m_j+n_j)$, $l_i=\sum_{j=1}^{i-1}(m_j+n_j)+n_i+1$,
$i=1,\ldots,N$.
Let $\eta$ be the configuration from the class of equivalence~$S$
such that $\eta(x)=1$ for $x\leq 0$ and $\eta(1)=0$.
  Define the configurations $\eta^\to_k$, $\eta^\gets_k$, 
$\eta^{+r}_k$,
$\eta^{+l}_k$, $\eta^{-r}_k$, $\eta^{-l}_k$ in the following way:
\begin{itemize}
\item $\eta^\to_k(x) = \eta(x)$ for $x\neq r_k, r_k+1$, $\eta^\to_k(r_k)=0$,
$\eta^\to_k(r_k+1)=1$, $k=0,\ldots,N$;
\item $\eta^\gets_k(x) = \eta(x)$ for $x\neq l_k, l_k-1$, 
$\eta^\gets_k(l_k)=0$,
$\eta^\gets_k(l_k-1)=1$, $k=1,\ldots,N$;
\item $\eta^{+r}_k(x) = \eta(x)$ for $x\neq r_k+1$, 
$\eta^{+r}_k(r_k+1)=1$, $k=0,\ldots,N$;
\item $\eta^{+l}_k(x) = \eta(x)$ for $x\neq l_k-1$,
$\eta^{+l}_k(l_k-1)=1$, $k=1,\ldots,N$;
\item $\eta^{-r}_k(x) = \eta(x)$ for $x\neq r_k$, 
$\eta^{-r}_k(r_k)=0$, $k=0,\ldots,N$;
\item $\eta^{-l}_k(x) = \eta(x)$ for $x\neq l_k$,
$\eta^{-l}_k(l_k)=0$, $k=1,\ldots,N$.
\end{itemize}
and $S^\to_k$, $S^\gets_k$, 
$S^{+r}_k$,
$S^{+l}_k$, $S^{-r}_k$, $S^{-l}_k$ are the corresponding
classes of equivalence.
Informally speaking,
\begin{itemize}
\item $S^\to_k$ is the configuration obtained from $S$ by moving the
  rightmost $1$ of the $k$-th $1$-block by $1$ unit to the right,
  $k=0,\ldots,N$;
\item $S^\gets_k$ is the configuration obtained from $S$ by moving the
  leftmost $1$ of the $k$-th $1$-block by $1$ unit to the left,
  $k=1,\ldots,N$;
\item $S^{+r}_k$ is the configuration obtained from $S$ by adding an
  extra $1$ to the right of the $k$-th $1$-block, $k=0,\ldots,N$;
\item $S^{+l}_k$ is the configuration obtained from $S$ by adding an
  extra $1$ to the left of the $k$-th $1$-block, $k=1,\ldots,N$;
\item $S^{-r}_k$ is the configuration obtained from $S$ by removing the
  rightmost $1$ from the $k$-th $1$-block, $k=0,\ldots,N$;
\item $S^{-l}_k$ is the configuration obtained from $S$ by removing the
  leftmost $1$ from the $k$-th $1$-block, $k=1,\ldots,N$.
\end{itemize}
Clearly, EP can transform $S$ to $S^\to_k$ or $S^\gets_k$, while using
VM we can get $S^{\pm r}_k$ or $S^{\pm l}_k$.

Denote also $R_i = \sum_{j=1}^i n_j$, $T_i = \sum_{j=i}^N m_j$, and let
$$
|S| = \sum_{j=1}^N (m_j+n_j) = R_N + T_1 
$$
stand for the length of ``nontrivial'' part of
configuration~$S$. Notational convention: $R_0 = T_{N+1} = 0$.

 We define
two functions $f_1,f_2: \tilde\cD \mapsto \R$, which will play
the crucial role in our arguments:
\begin{eqnarray*}
 f_1(S) & = & \frac12 \bigg(\sum_{k: S(k)=1}\Big(\sum_{m<k}
 \Ind{S(m)=0}\Big)+
\sum_{k: S(k)=0}\Big(\sum_{m>k}\Ind{S(m)=1}\Big)\bigg) \\
 &=& \frac12 \Big(\sum_{i=1}^N m_i R_i + 
            \sum_{i=1}^N n_iT_i \Big) \\
 &=& \sum_{i=1}^N m_i R_i = \sum_{i=1}^N n_iT_i,
\end{eqnarray*}
and
\begin{eqnarray*}
 f_2(S) & = & \frac12 \bigg(\sum_{k: S(k)=1}\Big(\sum_{m<k}
 \Ind{S(m)=0}\Big)^2 +
\sum_{k: S(k)=0}\Big(\sum_{m>k}\Ind{S(m)=1}\Big)^2\bigg) \\
 &=& \frac12 \bigg(\sum_{i=1}^N m_i R_i^2 + 
            \sum_{i=1}^N n_iT_i^2 \bigg),
\end{eqnarray*}
for all $S\in \tilde\cD$.

Before going further, let us make some remarks about $f_1$, $f_2$.
The value~$f_1(S)$ is equal exactly to the number of
nearest-neighbor transpositions needed to pass from~$S$ to~$\cD_0$,
that is, $f_1(S)$ is in some sense the ``distance'' from~$S$
to the trivial configuration. Unfortunately, as we will see later,
the function~$f_1$ does not ``work'' well for some configurations~$S$
(namely, for~$S$ such that $N(S)$ is small with respect to~$|S|$).
The function~$f_2$ is the result of our attempts to modify~$f_1$
in order to eliminate this disadvantage; we cannot give any
intuitive meaning of~$f_2(S)$.

Let us obtain some relations between $|S|$, $f_1(S)$ and
$f_2(S)$.

\begin{lmm}
\label{relations}
 For any $S\in\cD$ the following holds:
\begin{itemize}
\item[{\rm i)}] $|S|/2 \leq f_1(S) \leq |S|^2/4$;
\item[{\rm ii)}] $|S|^2/4 \leq f_2(S) \leq |S|^3/8$;
\item[{\rm iii)}] $f_1(S) \leq \big(f_2(S)\big)^{3/4}$.
\end{itemize}
\end{lmm}

\noindent
{\it Proof.} The proof of i)--ii) is simple. 
We have
$$ f_1(S)=\frac12 \Big(\sum_{i=1}^N m_i R_i + 
            \sum_{i=1}^N n_iT_i \Big) \geq \frac12 (R_N + T_1) =
   \frac{|S|}{2}, $$
$$ f_1(S)=\sum_{i=1}^N m_i R_i \leq R_N \sum_{i=1}^N m_i =
     R_NT_1 \leq \frac{(R_N+T_1)^2}{4}= \frac{|S|^2}{4}, $$
and, analogously,
$$ f_2(S) \geq \frac12 (R_N^2 + T_1^2) \geq
   \frac14(R_N+T_1)^2   =
   \frac{|S|^2}{4}, $$
$$ f_2(S)\leq \frac12 \Big(R_N^2\sum_{i=1}^N m_i + 
            T_1^2\sum_{i=1}^N n_i \Big) =
   \frac12 R_NT_1(R_N + T_1) \leq \frac{|S|^3}{8}.$$

Let us prove iii). We shall make use of the following
simple consequence of the Jensen inequality: if we have~$n$
positive numbers $\gamma_1,\ldots,\gamma_n$ such that
$\sum_{i=1}^n \gamma_i = 1$, then for any $x_1,\ldots,x_n$
\be
\label{jensen}
\gamma_1x_1+\cdots+\gamma_nx_n \leq 
    (\gamma_1x_1^2+\cdots+\gamma_nx_n^2)^{1/2}.
\ee
Denote $\alpha_i = m_i/|S|$, $\beta_i = n_i/|S|$, so
$\sum_{i=1}^N(\alpha_i+\beta_i) = 1$. 
Using~(\ref{jensen}) and ii), we get
\begin{eqnarray*}
f_1(S) & = & \frac12 \sum_{i=1}^N (m_iR_i+n_iT_i) = 
        \frac{|S|}{2} \sum_{i=1}^N (\alpha_iR_i+\beta_iT_i)\\
 &\leq & \frac{|S|}{2} \Big(\sum_{i=1}^N (\alpha_iR_i^2+\beta_iT_i^2)
   \Big)^{1/2} = \frac{\sqrt{|S|}}{\sqrt{2}} \big(f_2(S)\big)^{1/2}\\
  &\leq & \frac{\sqrt{2}\big(f_2(S)\big)^{1/4}}{\sqrt{2}}
   \big(f_2(S)\big)^{1/2} = \big(f_2(S)\big)^{3/4},
\end{eqnarray*}
thus completing the proof of Lemma~\ref{relations}.
\qed

\medskip

As usual, symbols $\P $ and $\E$ stand for probability
and expectation. When using them may look ambiguous, we
use symbol $\Ee{p}$ ($\Pe{p}$) to denote expectation (probability)
w.r.t.\ EP($p$), 
$\Ev$ ($\Pv$) stands for expectation (probability)
w.r.t.\ VM,
$\Eh{\beta}{p}$ ($\Ph{\beta}{p}$) denotes
 expectation (probability)
w.r.t.\ HP($\beta$, $p$).

\section{Exclusion process}
\label{excl}
In this section we shall study the EP 
using the method of Lyapunov functions.

\begin{theo}
\label{ex_erg}
 If $p>q$, then the exclusion process
is ergodic.
\end{theo}

\noindent{\it Proof.} 
As we noticed before, EP can transform a configuration $S$
only either to $S_k^\to$ or to $S_k^\gets$, where the notations
$S_k^\to$ and $S_k^\gets$ have been introduced in Section~\ref{notation}.
Then, it is elementary to get that
\begin{eqnarray}
 f_2(S^\to_k) - f_2(S)& = & \frac12 \big( (R_k+1)^2-R_k^2 +
 (T_{k+1}+1)^2
         - T_{k+1}^2 \big)\nonumber\\
& =& 1+R_k+T_{k+1} \label{eq_to}
\end{eqnarray}
and
\begin{eqnarray}
 f_2(S^\gets_k) - f_2(S)& = & \frac12 \big((R_k-1)^2-R_k^2 + 
(T_{k}-1)^2
         - T_{k}^2 \big)\nonumber\\
& =& 1-R_k-T_{k} \label{eq_gets}
\end{eqnarray}

Combining (\ref{eq_to}) and (\ref{eq_gets}), we have that

\be
\label{xxx}
\E(f_2(\xi_{t+1})-f_2(\xi_t) \mid \xi_t = S)
={N+q\over 2N+1}-{p-q\over 2N+1}
\sum_{i=1}^N (R_i+T_i). 
\ee
Since $R_N+T_1=|S|$, $R_i \geq i$ and $T_i \geq N-i+1$,
it is straightforward to get that
$$ \sum_{i=1}^N (R_i+T_i) \geq \max\{|S|,N(N+1)\}.$$
Using this fact, we get from~(\ref{xxx}) that for any~$\eps>0$
\be
\label{final_e}
  \E(f_2(\xi_{t+1})-f_2(\xi_t) \mid \xi_t = S) < -\eps  
\ee
for all but finitely many~$S$. So, 
by Theorem~\ref{cr_erg}, EP($p$) is ergodic when
$p>1/2$.
\qed

\begin{theo}
\label{ex_tr}
When $p\leq q$ the exclusion process
is transient.
\end{theo}

\noindent{\it Proof.} 
First we consider the case $p<q$.

With $S^\to_k$ and $S^\gets_k$ being as defined above, we have that
\be
\label{eq_fgets}
f_1(S^\gets_k)- f_1(S)=-1, 
\ee
and
\be
\label{eq_fto}
f_1(S^\to_k)- f_1(S)=1, 
\ee
so that for some $\eps=\eps(p,q)>0$
\be
\label{f_1ex}
\Ee{p}(f_1(\xi_{t+1})-f_1(\xi_t) \mid \xi_t = S)
 =  {{N(q-p)}\over{2N+1}}+{{q}\over{2N+1}}\geq
\eps
\ee
and, clearly,
$ |f_1(\xi_{t+1})-f_1(\xi_t)|\leq 1$
almost surely.
Then by Theorem~\ref{cr_tr2}, the process $\xi_t$ is transient.

Let us turn now to the case $p=q=1/2$.

Using the function $f_1(S)$ defined above and (\ref{eq_fgets}),
(\ref{eq_fto}), we have
that 
\be
\label{=frac}
\Ee{1/2}(f_1(\xi_{t+1})-f_1(\xi_t) \mid \xi_t = S)
= \frac{1}{2(2N+1)}
\ee
so Theorem~\ref{cr_tr1} does not apply.
Therefore, we need a different approach.

We fix an arbitrary $\alpha>0$ and define the function
$\psi : \cD\setminus\cD_0 \mapsto \R$ by
$$
\psi(S):=\bigl(f_1(S)\bigr)^{-\alpha}.
$$
Note that the definition is correct because $f_1(S)>0$ for $S\notin\cD_0$.
(Actually, for the need of Theorem~\ref{ex_tr} it is sufficient
 to take $\alpha=1$, but, since we
will need analogous calculations later in this paper,
at this point we prefer to do
the calculations for arbitrary~$\alpha>0$.)
To study the properties of the process $\psi(\xi_t)$, we
need the following lemma.

\begin{lmm}
\label{f>N}
For any $C>0$ the set
\be
\label{fin_set}
A_C = \{S : f_1(S)<CN(S)\}
\ee
is finite.
\end{lmm}

\noindent
{\it Proof.} Clearly, $R_i\geq i$ and $m_i\geq 1$, 
so $f_1(S)\geq N(S)(N(S)+1)/2$.
Thus, for a configuration~$S$ to belong to $A_C$,
it is necessary that the number of $1$-blocks
 be less than $2C-1$, so
 $A_C$ is a subset of
$$ \{S : f_1(S)<C(2C-1)\},$$
which is obviously finite.
\qed

It follows from~(\ref{eq_fgets}) and~(\ref{eq_fto}) that
\be
\label{=1}
\Ee{1/2}((f_1(\xi_{t+1})-f_1(\xi_t))^2 \mid \xi_t = S) = \frac12.
\ee
  By elementary calculations, we get that for any $\alpha>0$ there exist two
positive numbers $C_1=C_1(\alpha)$, $C_2=C_2(\alpha)$ such that
\be
\label{calc}
(x+1)^{-\alpha}-1 \leq  - \alpha x + C_1 x^2,
\ee
when $|x|<C_2$.

Using (\ref{=frac}), (\ref{=1}), (\ref{calc}) and Lemma~\ref{f>N} we get
\begin{eqnarray}
\lefteqn{\Ee{1/2}(\psi(\xi_{t+1})-\psi(\xi_t)\mid \xi_t = S)}
\nonumber\\
&=& f_1^{-\alpha}(S) \Ee{1/2}\bigg(\Big(\frac{f_1(\xi_{t+1})}
{f_1(\xi_t)}\Big)^{-\alpha}
      - 1 \mid \xi_t=S\bigg)\nonumber \\
&=& f_1^{-\alpha}(S) \Ee{1/2}\bigg(\Big(\frac{f_1(\xi_{t+1})-
f_1(\xi_t)}{f_1(\xi_t)} 
        +1  \Big)^{-\alpha} -1 \mid \xi_t=S\bigg)\nonumber \\
&\leq & f_1^{-\alpha}(S) \Big(-\frac{\alpha}{f_1(S)}\cdot \frac{1}{2(2N+1)} + 
               \frac{C_1}{2f_1^2(S)}\Big)\nonumber\\
&=& f_1^{-\alpha-2}(S) \Big(-\frac{\alpha f_1(S)}{2(2N+1)} + 
               \frac{C_1}{2}\Big) < 0 \label{taylor}
\end{eqnarray}
on $\big\{S: f_1(S) > \max\{1/C_2,C_1(2N(S)+1)/\alpha\}\big\}$, and hence
for all but finitely many~$S$. Applying Theorem~\ref{cr_tr2},
we finish the proof of Theorem~\ref{ex_tr}.
\qed

\section{Voter model}
\label{voter}
The subject of this section is the discrete time
voter model.
For the process starting from a configuration $S$ denote by
$\tau(S)$ the moment of hitting the set~$\cD_0$.
The main result of this section is the following
\begin{theo}
\label{t_vot}
The discrete time voter model is positive recurrent. Moreover,
for any initial configuration $S_0$ and any $\eps>0$
\be
\label{mom<infty}
 \E (\tau(S_0))^{3/2-\eps} < \infty
\ee  
and
\be
\label{mom=infty}
 \E (\tau(S_0))^{3/2+\eps} = \infty.
\ee
\end{theo}

\noindent
{\it Proof.} Since positive recurrence means just the existence
of $\E\tau(S_0)$, we shall turn directly to the proof
of~(\ref{mom<infty}). The idea is to apply Theorem~\ref{cr_mom}
to the process $f_2^\alpha(\xi_t)$ for some $\alpha<1$.

First, we need the following important fact

\begin{lmm}
\label{Ef_2=0}
We have
\be
\label{eq_Ef_2=0}
\Ev (f_2(\xi_{t+1})-f_2(\xi_t) \mid \xi_t=S) = 0
\ee
for any $S\in \cD$.
\end{lmm}

\noindent
{\it Proof.} If $S\in \cD_0$, then~(\ref{eq_Ef_2=0}) is
trivial. For $S\notin \cD_0$ a direct computation
gives
\begin{eqnarray}
 f_2(S_k^{+r})-f_2(S) &=& \frac12 (R_k+T_{k+1}+R_k^2-T_{k+1}^2) -
     \sum_{i=k+1}^N m_iR_i +\sum_{i=1}^k n_iT_i,\phantom{***} 
          \label{eee1}\\
  f_2(S_k^{-r})-f_2(S) &=& \frac12 (R_k+T_{k+1}-R_k^2+T_{k+1}^2) +
     \sum_{i=k+1}^N m_iR_i -\sum_{i=1}^k n_iT_i \phantom{***}
          \label{eee2}
\end{eqnarray}
for $k=0,\ldots,N$, and
\begin{eqnarray}
f_2(S_k^{+l})-f_2(S) &=& \frac12 (-R_k-T_{k}+R_k^2-T_{k}^2) -
     \sum_{i=k}^N m_iR_i +\sum_{i=1}^k n_iT_i,\phantom{***}
               \label{eee3}\\
f_2(S_k^{-l})-f_2(S) &=& \frac12 (-R_k-T_{k}-R_k^2+T_{k}^2) +
     \sum_{i=k}^N m_iR_i -\sum_{i=1}^k n_iT_i .\phantom{***}
         \label{eee4}
\end{eqnarray}
for $k=1,\ldots,N$.

Taking summation in (\ref{eee1})--(\ref{eee4}) one gets~$0$,
thus finishing the proof of Lemma~\ref{Ef_2=0}. \qed

\medskip

Then, from~(\ref{eee2}) we note that 
\be
\label{vvv1}
|f_2(S_0^{-r})-f_2(S)| \geq \frac{T_1^2}{2}
\ee
and from~(\ref{eee1})
\be
\label{vvv2}
|f_2(S_N^{+r})-f_2(S)| \geq \frac{R_N^2}{2}.
\ee
These two inequalities give us that
there exist a constant~$C>0$ such that
\be
\label{2mom_f_2'}
\Ev ((f_2(\xi_{t+1})-f_2(\xi_t))^2 \mid \xi_t=S) \geq \frac{C|S|^4}{N}
\ee
for all $S$.
Now, a very important observation is that the VM does not
increase the number of blocks $N_t=N(\xi_t)$. So we have for all~$S$
\be
\label{2mom_f_2}
\Ev ((f_2(\xi_{t+1})-f_2(\xi_t))^2 \mid \xi_t=S) \geq C_0|S|^4
\ee
with $C_0 = C_0(S_0) = C/N(S_0)$. 

Elementary calculus gives us that for $0<\alpha <1$ and for
$|x|\leq 1$ there exists a positive constant $C_1$ such that
\be
\label{aaaa}
 (x+1)^\alpha - 1 \leq \alpha x - C_1 x^2. 
\ee

Using now considerations analogous to~(\ref{taylor})
and applying~(\ref{aaaa}), Lemma~\ref{Ef_2=0} and~(\ref{2mom_f_2})
we get
\be
\label{ocenka}
 \Ev ((f_2(\xi_{t+1}))^\alpha-(f_2(\xi_t))^\alpha \mid \xi_t=S) \leq
     - C_0 C_1 (f_2(S))^{\alpha-2} |S|^4. 
\ee
Applying Lemma~\ref{relations}, part ii), to the last inequality 
we get
$$ \Ev ((f_2(\xi_{t+1}))^\alpha-(f_2(\xi_t))^\alpha \mid \xi_t=S) \leq
     - 16C_0 C_1(f_2(S))^{\alpha-2/3} . $$
We apply Theorem~\ref{cr_mom} to the process $X_t = (f_2(\xi_t))^{1/3}$
taking~$\alpha$ to be close to~$1$
to finish the proof of~(\ref{mom<infty}).

\medskip

Let us turn now to the proof of~(\ref{mom=infty}).
We let the process start from configuration~$S_0$
such that $N(S_0)=1$. Since this
configuration is reachable from any other configuration,
it is sufficient to prove~(\ref{mom=infty}) for
this~$S_0$. As it was mentioned before, the voter
model does not make the number of blocks~$N$ increase,
so the process can be represented as $\xi_t=(n_t,m_t)$,
which clearly is a random walk in $\Z_+^2$, and we are
interested in the moment of hitting the boundary.
Note that the transition probabilities of this random
walk can be described like this: from the state $(n,m)$
the transition can occur to the states $(n+1,m)$, $(n-1,m)$,
$(n,m+1)$, $(n,m-1)$, $(n+1,m-1)$ and $(n-1,m+1)$ with
probabilities $1/6$.

Denote by $\tau_{n,m}$ the moment of hitting $\cD_0$
(i.e.\ the boundary) provided that the starting point was 
$(n,m)$. To proceed, we need the following

\begin{lmm}
\label{lamperti}
There exist two positive constants $\delta$, $C$, such
that for any $n$, $m$
\be
\label{n^2}
\P\{\tau_{n,m}>\delta n^2\} \geq \frac{Cm}{m+n}
\ee
and
\be
\label{m^2}
\P\{\tau_{n,m}>\delta m^2\} \geq \frac{C n}{m+n} .
\ee
\end{lmm}

\begin{rem}
It can be shown that Lemma~\ref{lamperti} holds for
any homogeneous random walk in $\Z^2_+$ with bounded
jumps and zero drift in the interior of $\Z^2_+$.
\end{rem}

\noindent
{\it Proof.} Without loss of generality we can suppose
that $n\leq m$. Then, to prove~(\ref{n^2}), we will
prove a stronger fact:
\be
\label{stron}
\P\{\tau_{n,m}>\delta n^2\} \geq C_0
\ee
for some $C_0$.  
In fact, it is a  classical result that a 
homogeneous random walk in $\Z^2_+$ with bounded
jumps and zero drift in the interior
with some uniformly positive probability cannot deviate
by the distance~$n$ from its initial position during
the time~$n^2$. To show how it can be proved formally,
we denote by $\rh{(n_1,m_1),(n_2,m_2)}$ the
Euclidean distance between the points $(n_1,m_1)$
and $(n_2,m_2)$. Let the process start from $(n,m)$, 
and  denote $Y_t = \rh{\xi_t, (n,m)}$. Then, it is
straightforward to get that the process~$Y_t$ 
satisfies the hypothesis of Lemma~2 from Aspandiiarov \ea
(1996), so
applying it, we finish the proof of~(\ref{n^2}).

 To prove~(\ref{m^2}), we need some additional notations.
Denote
\begin{eqnarray*}
W^i(m) & = & \{(n',m') : \rh{(n',m'),(m,m)}\leq m/2 + \sqrt{2}\},\\
V^i(m) & = & \{(n',m') : m/2 <\rh{(n',m'),(m,m)}\leq m/2 + \sqrt{2}\},\\
W^e(m) & = & \{(n',m') : m/2 + \sqrt{2}<\rh{(n',m'),(m,m)}\leq m\},\\
V^e(m) & = & \{(n',m') : m -\sqrt{2} < \rh{(n',m'),(m,m)}\leq m \}.
\end{eqnarray*}
Clearly, the set $V^i(m)$ is the boundary of $W^i(m)$, and
the set  $V^e(m)$  is the external boundary of $W^e(m)$.

 We consider the two possible cases:
\begin{itemize}
\item[{\rm a)}] $(n,m) \in W^i(m)$,
\item[{\rm b)}] $(n,m) \in W^e(m)$.
\end{itemize}

 Case a): first, we denote $Y_t=\rh{\xi_t,(m,m)}$. Then,
we apply Lemma~2 from Aspandiiarov \ea
(1996) to get that 
$\P\{\tau_{n,m}>\delta n^2\} \geq C_1$ for some~$C_1$, and
thus~(\ref{m^2}).

 Case b): we keep the notation $Y_t$ from the previous
paragraph. Denote by $p_{n,m}$ the probability of hitting
the set $V^i(m)$  before the set $V^e(m)$,  provided
that the starting point is $(n,m)$. Our goal is to 
estimate this probability from below.

 For $C>0$ consider the process $Z^C_t$ , $t=0,1,2,\ldots$,
defined in the following way:
$$Z^C_t = \exp\Big\{C\Big(1-\frac{Y_t}{m}\Big)\Big\} = 
   \exp\Big\{C\Big(1-\frac{\rh{\xi_t,(m,m)}}{m}\Big)\Big\},$$
$Z^C_0 = \exp\{Cn/m\}$. One can prove the 
following technical fact: there exists a constant~$C$
(not depending on~$m$) such that 
\be
\label{mart11}
\E(Z^C_{t+1} - Z^C_t \mid \xi_t = (n',m')) \geq 0
\ee
for any point $(n',m')\in W^e(m)$ and if~$m$ is large enough. 
Indeed, using the fact that
there exist two positive constants $C_{1,2}$ such that
$$ e^{-x} - 1 \geq -x + C_1 x^2$$ 
on $|x|<C_2$, we write
\begin{eqnarray*}
\lefteqn{\E(Z^C_{t+1} - Z^C_t \mid \xi_t = (n',m'))}\\
&=& \exp\Big\{C\Big(1-\frac{|n'-m|}{m}\Big)\Big\}
\E\Big(\exp\Big\{-\frac{C}{m}(Y_{t+1}-Y_t)\Big\}-1
          \mid \xi_t = (n',m') \Big)\\
&\geq & \frac{C}{m} \exp\Big\{C\Big(1-\frac{|n'-m|}{m}\Big)\Big\}
 \E\Big(-(Y_{t+1}-Y_t) + \frac{C_1C}{m} (Y_{t+1}-Y_t)^2
 \mid \xi_t = (n',m')\Big).
\end{eqnarray*}
Then, using properties of the process~$Y_t$, one can
complete the proof of~(\ref{mart11}).

Now, to estimate $p_{n,m}$, we make the sets $V^i(m)$
and $V^e(m)$ absorbing. Using that our random walk cannot
overpass these sets, the 
process $Z^C_t$ converge as $t\to\infty$ to $Z^C_\infty$, so
\begin{eqnarray*}
\E Z^C_\infty &\geq& p_{n,m}\exp\Big\{\frac{C}{2}\Big\} + 
                  (1-p_{n,m}) \\
      & \geq & \E Z^C_0 = \exp\Big\{\frac{Cn}{m}\Big\},
\end{eqnarray*}
and thus
\begin{eqnarray}
p_{n,m} &\geq & \frac{\exp\{Cn/m\}-1}{\exp\{C/2\}-1} \nonumber \\
&\geq &\frac{C}{\exp\{C/2\}-1} \cdot \frac{n}{m} \geq
    \frac{2C}{\exp\{C/2\}-1} \cdot \frac{n}{m+n}. \label{occ}
\end{eqnarray}
So, starting from the point~$(n,m)$, with probability
at least~(\ref{occ}) the random walk hits the set $V^i(m)$.
 Then, from the case a) it follows that with
uniformly positive probability it will take at least~$\delta m$
steps to reach the external boundary~$V^e(m)$, so we complete the proof
of~(\ref{m^2}) and thus, of Lemma~\ref{lamperti}.
\qed

\medskip

Now, supposing that~(\ref{mom=infty}) does not hold, 
 we have (denoting $\tau := \tau(S_0)$ and $a\wedge b := \min\{a,b\}$)
\begin{eqnarray}
\E \tau^{3/2+\eps}
&\geq &  \E( \tau^{3/2+\eps}\Ind{\tau \geq t}) =
\E ( (t+\tau_{\xi_t})^{3/2+\eps}\Ind{\xi_s \notin\cD_0 
     \mbox{ \scriptsize for all }s\leq t}
     )\nonumber\\
&\geq & \frac12 \E \Big( (t+\delta n^2_t)^{3/2+\eps}\frac{Cm_t}{m_t+n_t}
\Ind{\xi_s\notin\cD_0 
 \mbox{ \scriptsize for all }s\leq t}
      \Big)\nonumber\\
& &{} +
\frac12 \E \Big( (t+\delta m^2_t)^{3/2+\eps}\frac{Cn_t}{m_t+n_t}
\Ind{\xi_s\notin\cD_0
 \mbox{ \scriptsize for all }s\leq t}
 \Big)\nonumber\\
&\geq &  \delta' C' \E\big(\bigl(
n^{2+\eps'}_t m_t+m^{2+\eps'}_tn_t\bigr)
\Ind{\xi_s\notin\cD_0 
\mbox{ \scriptsize for all }s\leq t}\big)\nonumber\\
&= &  \delta' C' \E\big(\bigl(
n^{2+\eps'}_{t\wedge\tau} m_{t\wedge\tau}+
  m^{2+\eps'}_{t\wedge\tau}n_{t\wedge\tau}\bigr)
   \big)\nonumber\\
&=& C^{\prime\prime} \E (f_2 (\xi_{t\wedge\tau}))^{1+\epss}.
\label{longstory}
\end{eqnarray}
for some constants $\delta'$, $\eps'$, $C'$,
 $\epss$ and $C^{\prime\prime}$.

  From (\ref{longstory})  we get
that the family $\{f_2 (\xi_{t\wedge\tau})\}$ is uniformly
integrable as $t\to\infty$, so
 $\E f_2(\xi_t) \to \E f_2(\xi_\tau) =  0$. But this obviously
contradicts to Lemma~\ref{Ef_2=0}.
\qed

\section{Hybrid process}
\label{hybrid}
As it was proved before, the EP($p$) is transient when $p\leq 1/2$,
and VM is ergodic. Now, what will happen if we combine them?
The following theorems give a (not complete) answer to
this question.

\begin{theo}
\label{hyb_erg}
 There exists $\beta_0<1$
 such that for any $p$ the process HP($\beta$, $p$) is
ergodic for all $\beta>\beta_0$.
\end{theo}

\begin{theo}
\label{hyb_symm}
For any $\beta>0$ and $p\geq 1/2$ 
the process HP($\beta$, $1/2$) is ergodic. 
\end{theo}

We also formulate the following plausible conjecture.
Its not completely rigorous proof will be presented in 
Section~\ref{hp_nerg}.
\begin{conj}
\label{cnj}
For any $p<1/2$ there exists $\beta_0=\beta_0(p)>0$ such that
the process HP($\beta$, $p$) is not ergodic for $\beta<\beta_0$.
\end{conj}

\subsection{Proof of Theorem~\protect\ref{hyb_erg}}
\label{proof_hyb_erg}
We use the notations introduced in Section~\ref{notation}.
Direct computations yield
\begin{eqnarray*}
f_1(S_k^{+l}) - f_1(S) &=& R_k - T_k -1, \\
f_1(S_k^{-l}) - f_1(S) &=& -R_k + T_k -1, \\
f_1(S_k^{+r}) - f_1(S) &=& R_k - T_{k+1} , \\
f_1(S_k^{+r}) - f_1(S) &=& - R_k + T_{k+1} ,
\end{eqnarray*}
so
\be
\label{f_1vot}
\Ev(f_1(\xi_{t+1}) - f_1(\xi_{t}) \mid \xi_t=S) = - \frac{N}{2N+1}.
\ee

Combining this with (\ref{f_1ex}), we get
that there exists a positive number $C=C(\beta)$
such that
\begin{eqnarray}
\Eh{\beta}{p}(f_1(\xi_{t+1}) - f_1(\xi_{t}) \mid \xi_t=S)& =&
 - \frac{1}{2N+1} \big(\beta N - (1-\beta)
     ((q-p) N + q)\big) \nonumber \\
 & < & - C(\beta) \label{hyb_beta_p}
\end{eqnarray}
for $\beta>2/3$. Applying Theorem~\ref{cr_erg},
we finish the proof. \qed

\subsection{Proof of Theorem~\protect\ref{hyb_symm}}
\label{proof_hyb_symm}

To prove the desired result, we are going to apply
 Theorem~\ref{cr_erg} to the function $\phi(S):=
  (f_2(S))^\alpha$ for some $\alpha<1$. 

 First, we prove the theorem for the case $p=1/2$.

 Inserting $p=q=1/2$ into~(\ref{xxx}), we obtain
for the step of EP($1/2$) 
\be
\label{ex_1_mom}
\Ee{1/2}(f_2(\xi_{t+1}) - f_2(\xi_{t}) \mid \xi_t=S) =
        \frac{1}{2}.
\ee
It is elementary to get that for $\alpha\in (0,1)$
\be
\label{bern} 
(x+1)^\alpha - 1 \leq \alpha x
\ee
for all $x\geq -1$.
  Using (\ref{ex_1_mom}), (\ref{bern}) and~(\ref{aaaa}),
we get
\be
\label{ocenka_ex}
\Ee{1/2}((f_2(\xi_{t+1}))^\alpha - (f_2(\xi_{t}))^\alpha \mid \xi_t=S)
   \leq \frac{\alpha (f_2(S))^{\alpha-1}}{2} . 
\ee

Now, let us make the necessary estimate for the step of VM.
Here we will need a bound which is more accurate 
than~(\ref{2mom_f_2'}):
\begin{lmm}
\label{l_16/5}
There exists $C'>0$ such that
\be
\label{16/5}
\Ev((f_2(\xi_{t+1} - f_2(\xi_{t}))^2 \mid \xi_t=S) \geq C' |S|^{16/5}
\ee
\end{lmm}

\noindent
{\it Proof.} To calculate exactly
 the left-hand side of~(\ref{16/5}), one has
to square (\ref{eee1})--(\ref{eee4}), sum them up and divide 
by~$4N+2$. But this calculation appears to be too difficult;
so we will only obtain a lower bound. Denote $\Delta_k = 
f_2(S_k^{+r}) - f_2(S)$, so
\be
\label{16/5**}
\Ev((f_2(\xi_{t+1} - f_2(\xi_{t}))^2 \mid \xi_t=S) \geq
    \frac{1}{4N+2}\sum_{i=1}^N\Delta_i^2.
\ee
By simple algebraic calculations, one gets from~(\ref{eee1}) that
\be
\label{16/5***}
\Delta_{i+1} - \Delta_i \geq N
\ee
for $i=0,\ldots, N-1$. From~(\ref{eee1}) one gets also that $\Delta_0<0$
and $\Delta_N>0$, so denote
 $L=\min\{k: \Delta_{k-1}<0, \Delta_k\geq 0\}$. Using~(\ref{16/5***}),
we get
\begin{eqnarray*}
\sum_{i=1}^N\Delta_i^2 &\geq & \sum_{i=0}^{L-1}(N(L-i-1))^2 +
   \sum_{i=L}^N (N(i-L))^2\\
 &\geq & N^2 \sum_{i=1}^{N/2}i^2 \geq C_1 N^5
\end{eqnarray*}
for some $C_1$, so by~(\ref{16/5**}) we get that
$$
\Ev((f_2(\xi_{t+1} - f_2(\xi_{t}))^2 \mid \xi_t=S) \geq C_2 N^4
$$
for some $C_2$. Combining this with~(\ref{2mom_f_2'}), we get
$$
\Ev((f_2(\xi_{t+1} - f_2(\xi_{t}))^2 \mid \xi_t=S) \geq 
       \max\Big\{C_2 N^4, \frac{C|S|^4}{N}\Big\} \geq 
          C^{4/5}C_2^{1/5}|S|^{16/5},
$$  
thus proving Lemma~\ref{l_16/5}.
\qed

\begin{rem}
The exponent~$16/5$ in Lemma~\ref{l_16/5} is the best possible;
to see this, one may take a configuration~$S$ with $n_1=m_N=N^{5/4}$
and $n_2=\cdots=n_N=m_1=\cdots=m_{N-1}=1$ and compute the
left-hand side of~(\ref{16/5}).
\end{rem}

We continue proving Theorem~\ref{hyb_symm}.
Analogously to~(\ref{taylor}), using~(\ref{aaaa}) together
with Lemmas~\ref{Ef_2=0} and~\ref{l_16/5}, we get
for some positive constant~$C_2$
\be
\label{ocenka_N}
 \Ev ((f_2(\xi_{t+1}))^\alpha-(f_2(\xi_t))^\alpha \mid \xi_t=S) \leq
     - C_2 (f_2(S))^{\alpha-2} |S|^{16/5}.
\ee
So, for HP($\beta$,$1/2$), combining~(\ref{ocenka_ex}) with~(\ref{ocenka_N})
and using that $|S|^{16/5} \geq 2^{16/5}(f_2(S))^{16/15}$ because of
part ii) of Lemma~\ref{relations}, we get for $14/15 < \alpha < 1$
\begin{eqnarray}
\lefteqn{\Eh{\beta}{1/2}(\phi(\xi_{t+1})-\phi(\xi_t)
                                              \mid \xi_t=S)}
  \nonumber \\
 &=& \Eh{\beta}{1/2}((f_2(\xi_{t+1}))^\alpha-(f_2(\xi_t))^\alpha 
                                              \mid \xi_t=S)
  \nonumber \\
 &\leq & (1-\beta) \frac{\alpha (f_2(S))^{\alpha-1}}{2}
           - \beta C_2 (f_2(S))^{\alpha-2} |S|^{16/5}
          \nonumber   \\
 &\leq & (1-\beta) \frac{\alpha (f_2(S))^{\alpha-1}}{2}
           - 2^{16/5}\beta C_2 (f_2(S))^{\alpha-14/15}\nonumber \\
 &=& -(f_2(S))^{\alpha - 14/15} \Big[2^{16/5}\beta C_2 - 
     \frac{(1-\beta)\alpha}{2}(f_2(S))^{-1/15}\Big] < -1
\label{ocenka_hyb}
\end{eqnarray}
for all but a finite number of~$S$'s (indeed, the expression in
the square brackets is of order of positive constant for
 all but finitely many~$S$, 
and the fact that $\alpha>14/15$ guarantees that the
absolute value of the left-hand side of~(\ref{ocenka_hyb}) 
is large enough for all but a finite number of~$S$'s).
 Applying Theorem~\ref{cr_erg},
we finish the proof of Theorem~\ref{hyb_symm} for $p=1/2$.

 Now, when $p>1/2$, using~(\ref{final_e}) and Lemma~\ref{Ef_2=0}
we get that for any $\eps>0$
$$
 \Eh{\beta}{p}(f_2(\xi_{t+1})-f_2(\xi_t) \mid \xi_t=S) =
   (1-\beta)\Ee{p}(f_2(\xi_{t+1})-f_2(\xi_t) \mid \xi_t=S) <
     -(1-\beta)\eps
$$
for all but finite number of~$S$, and we apply Theorem~\ref{cr_erg}
again.
\qed

\begin{rem}
Using the technique of
Section~\ref{voter}, it is possible to get that there exists
some~$p_0=p_0(\beta)>1$ such that for the process 
 HP($\beta$, $1/2$) we have that $\E(\tau(S_0))^p < \infty$
for all $p<p_0$. By using the technique of Menshikov and Popov (1995),
one can get polynomial bounds on the decay of the stationary measure.
\end{rem}

\subsection{Nonergodicity}
\label{hp_nerg}
Here we will present an argument in the favor of the validity of
Conjecture~\ref{cnj}. 

 We rewrite~(\ref{hyb_beta_p}) as
\begin{eqnarray}
\Eh{\beta}{p}(f_1(\xi_{t+1}) - f_1(\xi_{t}) \mid \xi_t=S)& =&
  \frac{1}{2N+1} \big(-\beta N + (1-\beta)
     ((q-p) N + q)\big) \nonumber \\
 & \geq & \frac{N}{2N+1}\big((1-\beta)(q-p)-\beta\big) > 0
\label{hyb_beta_p>0}
\end{eqnarray}
when $\beta < \frac{q-p}{2q}$. Unfortunately,
because the VM does not have the property~(\ref{eq_tr2'}),
we cannot apply Theorem~\ref{cr_tr2}. Moreover, it is
still very unclear to us, if the process is transient
in this case. So instead we shall explain why we believe
it is not ergodic.
 
 We need the following three lemmas:

\begin{lmm}
\label{lll1}
 Let $\xi_t$, $t=0,1,2,\ldots$ be a
Markov chain on a countable state space~$X$,  let $0\in X$
be an absorbing state, and
define $ \tau :=  \min\{t: \xi_t=0 \}$ to be the hitting time of~$0$.
Suppose that for any starting point~$x$ we have $f(x):=\E_x\tau < \infty$.
 Then
\be
\label{L_1}
\E f(\xi_t) \to 0,
\ee
as $t\to\infty$.
\end{lmm}

\noindent
{\it Proof.} Let $x_0$ be the starting position of the
Markov chain. It is straightforward to get
\be
\label{fost}
\E (f(\xi_{t+1})-f(\xi_t) \mid \xi_t=x) = -\Ind{x\neq 0},
\ee
so, taking expectation in~(\ref{fost}), we get
\be
\label{fost'}
\E_{x_0} f(\xi_{t+1})-\E_{x_0} f(\xi_t) = -\P_{x_0}\{\tau > t\}.
\ee
 Taking summation in~(\ref{fost'}), we obtain
$$\E_{x_0} f(\xi_{t+1}) = f(x_0) - \sum_{i=0}^t\P_{x_0}\{\tau > i\} \to 0$$
as $t\to\infty$, thus completing the proof of Lemma~\ref{lll1}. 
\qed

\begin{lmm}
\label{lll2}
 Let $\xi_t$, $t=0,1,2,\ldots$ be a
Markov chain on a countable state space~$X$, and let $0\in X$
be an absorbing state. Let $x_0$ be the starting position of the
Markov chain, $\tau$ be the moment of hitting~$0$,
and suppose that $\E_{x_0}\tau<\infty$ for all~$x_0$.
 Let $f(x)$ be some nonnegative
function on~$X$ such that for some constant~$K$
\be
\label{men}
\E (f(\xi_{t+1})-f(\xi_t) \mid \xi_t=x) \leq  K
\ee
 for all $x\neq 0$. Then there exists a
constant~$M$ such that $\E f(\xi_t) < M$ for all~$t$.
\end{lmm}

\noindent
{\it Proof.} 
The proof is analogous to that of Lemma~\ref{lll1}:
first, we rewrite~(\ref{men}) as
\be
\label{men'}
\E (f(\xi_{t+1})-f(\xi_t) \mid \xi_t=x) \leq  K\Ind{\tau>t}
\ee
for all $x$. So,
\be
\label{men''}
\E_{x_0} f(\xi_{t+1})-\E_{x_0} f(\xi_t) \leq K\P_{x_0}\{\tau > t\},
\ee
and, taking summation in~(\ref{men''}), we get
$$\E_{x_0} f(\xi_{t+1}) \leq f(x_0) + K\sum_{i=0}^t\P_{x_0}\{\tau > i\} \leq
          f(x_0) + K\E_{x_0}\tau.    $$
Denoting $M:=f(x_0) + K\E_{x_0}\tau$, we finish the
proof of Lemma~\ref{lll2}.
\qed

\medskip

Analogously to Lemmas~\ref{lll1} and~\ref{lll2}, we can
prove the following lemma (which is, in fact, an adaptation
of Lemma~2.2 from Menshikov and Popov (1995) to our situation)

\begin{lmm}
\label{lll3}
 Let $\xi_t$, $t=0,1,2,\ldots$ be a
Markov chain on a countable state space~$X$, let  $0\in X$
be an absorbing state, $x_0$ be the starting point,
and $\tau$ be the moment of hitting~$0$,
and suppose that $\E_{x_0}\tau<\infty$ for all~$x_0$.
 Let $f(x)$ be some nonnegative
function on~$X$ such that $\E_{x_0} f(\xi_t)\to 0$ as $t\to\infty$,
and for some positive constant~$K$
\be
\label{menn}
\E (f(\xi_{t+1})-f(\xi_t) \mid \xi_t=x) \geq - K
\ee
 for all $x$. Then  $\E_{x_0} \tau \geq f(x_0)/K$.
\end{lmm}

 So, let us take a hybrid process HP($\beta$, $p$) which
satisfies~(\ref{hyb_beta_p>0}). We suppose that it is 
ergodic and try to get a contradiction. Define $\phi_{\beta,p}(S)$
to be the mean hitting time of~$\cD_0$ starting from~$S$, i.e.
$\phi_{\beta,p}(S) = \Eh{\beta}{p}\tau(S)$. We will prove
the following

\begin{lmm}
\label{<Cf_1}
For any $p>1/2$, $\beta$, there exists a positive constant
$C=C(\beta,p)$ such that 
\be
\label{cnjj}
\phi_{\beta,p}(S) \geq C f_1(S).
\ee
\end{lmm}

\noindent
{\it Proof.}   
  From (\ref{xxx}) and Lemma~\ref{Ef_2=0} we get that
for $p>1/2$ and any $\beta$ 
\be
\label{hhhhh}
\Eh{\beta}{p} (f_2(\xi_{t+1})-f_2(\xi_t)\mid \xi_t=S) \leq \frac 12.
\ee
 Applying Lemma~\ref{lll2}, we get that there exists a constant~$M$
such that $\Eh{\beta}{p} f_2(\xi_t) < M$ for all~$t$.
Using Lemma~\ref{relations} iii), we see that
$\Eh{\beta}{p} f_1(\xi_t) \to 0$ as $t\to\infty$. Applying
Lemma~\ref{lll3}, we complete the proof of Lemma~\ref{<Cf_1}.
\qed

\begin{conj}
\label{????}
Lemma~\ref{<Cf_1} holds for any $p$.
\end{conj}

We failed to prove the above conjecture. Intuitively,
$\phi_{\beta,p}(S)$ grows when~$p$ decrease and the monotonicity
argument might be applicable to prove this fact.
For the pure exclusion process, this argument follows
from the basic coupling (see Liggett (1985, Section~VIII.2)).
When the voter model is added, this coupling does not work.

 Now, if we suppose this to be true, the rest of the proof
is straightforward. If the process HP($\beta$, $p$) is ergodic,
then the function $\phi_{\beta,p}(S)$ is well defined, so,
by Lemma~\ref{lll1}, $\Eh{\beta}{p}\phi_{\beta,p}(\xi_t) \to 0$
 as $t\to\infty$.
Thus, using~(\ref{cnjj}), we get that $\Eh{\beta}{p}f_1(\xi_t)\to 0$.
But this obviously contradicts to~(\ref{hyb_beta_p>0}).
\qed

\section{Continuous time}\label{continuous}

In this section we show how Theorem~\ref{pp} follows from the theorems
proved in the last three sections.

Observe first that the transience is a property of the skeleton of a
Markov process. In our case the skeleton is the process
$\xi_n=\eta_{\tau_n}$, as defined in Section~\ref{discrete}. Notice
that according to this definition, $\xi_{n+1}$ may be the same as
$\xi_n$; this deviates a bit from the usual notion of skeleton. In our
version of skeleton the exit time of a configuration is a geometric
random variable with parameter bigger than $\beta +
(1-\beta)\min\{p,q\}$. This implies that the skeleton can not get
stacked. Hence Theorem~\ref{ex_tr} which states the transience for the
discrete-time exclusion process with $p\le 1/2$ implies the same for
the continuous-time process. This shows the transient part of item~1
of Theorem~\ref{pp}.

To prove that the ergodicity for the discrete-time process implies
the ergodicity for the continuous-time one, let $\eta\in\cD_0$, let $S$ be the
class of equivalence of $\eta$ and write
\be
  \label{p5}
  \tau_c(\eta)\,=\, \sum_{n=1}^{\tau(S)} \,(\tau_n-\tau_{n-1})
\ee
where we recall that $\tau_c(\eta)$ and $\tau(S)$ are the hitting
times of $\cD_0$ for the continuous and discrete time processes
starting from $\eta$ and $S$
respectively and $\tau_n$ is the instant of the $n$-th attempted jump
of the continuous process $\eta_t$, as defined in~(\ref{h2}). Given
the past up to $\tau_n$, $\tau_{n+1}-\tau_n$ is an exponential random
variable with rate bigger than 1 ---the worst case, the configurations
belonging to $\cD_0$. Hence, $\tau_{n+1}-\tau_n$ is stochastically
bounded below by an exponential random variable of rate 1
independent of everything. This implies that
\be
  \label{p5'}
  \E\tau_c(\eta)\,\le\, \,\E\tau(\eta).
\ee
Since the ergodicity is equivalent to the finiteness of the expected return
time for any given configuration then the  ergodicity for the discrete-time
process implies the same for the continuous-time one. With this
argument Theorem~\ref{ex_erg} implies the ergodic part of item~1 of
Theorem~\ref{pp} and Theorems~\ref{hyb_erg} and~\ref{hyb_symm} imply
item~2 of Theorem~\ref{pp}.

The argument above implies a stronger statement for the pure voter
model. Let $\beta=0$ and observe that for any $\eta \notin \cD_0$,
there are at least three discrepancies. Hence, for the voter model,
\be
  \label{p8}
  \tau_c(\eta)\,\le\, \sum_{n=1}^{\tau(S)} (\tau'_n-\tau'_{n-1}).
\ee
where $(\tau'_n-\tau'_{n-1})$ are independent exponentially distributed with
parameter $3$ and independent of $\tau(S)$.  This together with 
Lemma~\ref{p10} below show that the first part of~(\ref{p2}) follows 
from~(\ref{mom<infty}).

We now show how to obtain the second part of~(\ref{p2}) 
from~(\ref{mom=infty}). Let $\cD_1$ be the set of configurations on $\cD$
having exactly three discrepancies (that is, $10$, $01$ and $10$).
This means that in the representation~(\ref{m_in_i}) the
configurations belonging to $\cD_1$ have $N=1$, that is only one
finite 0-block and one finite 1-block.  The transition rate for
configurations in $\cD_1$ in the voter model is exactly 3. If the
process is in $\cD_1$ then it can only either stay in $\cD_1$ or 
jump to $\cD_0$.
Hence, for $\eta\in\cD_1$,
\be
  \label{p6}
 \tau_c(\eta)\,=\, \sum_{n=1}^{\tau(S)} (\tau'_n-\tau'_{n-1}).
\ee
Lemma~\ref{p10} below shows that~(\ref{mom=infty}) implies that the left
hand side of~(\ref{p6}) is infinite for any $\eta\in\cD_1$. As argued
before, any configuration in $\cD_1$ is reachable from any other,
hence the same is valid for any $\eta\in\cD\setminus\cD_0$. This shows
the second part of~(\ref{p2}).

\begin{lmm}
\label{p10}
Let $\tau$ be a positive integer random variable and $\tau_i$ be 
non\-ne\-ga\-tive independent
 random variables with the exponential distribution and
independent of~$\tau$. Then for any $p> 0$
 \be
   \label{p11}
   \E\Bigl(\sum_{n=1}^{\tau} \tau_n\Bigr)^p < \infty 
     \;\;\;\hbox{ if and only
   if }\;\;\; \E \tau^p <\infty.
 \ee
\end{lmm}

\noindent
{\it Proof.} By independence,
 \be
   \label{p12}
 \E\Bigl(\sum_{n=1}^{\tau} \tau_n\Bigr)^p =
   \sum_n \E\Bigl(\sum_{i=1}^{n} \tau_i \Bigr)^p \P(\tau=n).
 \ee
But 
\be
  \E\Bigl(\sum_{i=1}^{n} \tau_i \Bigr)^p \,=\, {\Gamma (n+p) \over
  \Gamma(n)} 
\ee
which is of the order of $n^p$.  \qed

\section{Final remarks}
Let us give several remarks with respect 
to extensions of our results to a nonnearest-neighbor
 case. For a nonnearest-neighbor voter model, we failed to
find an analogue of Lemma~\ref{Ef_2=0}. 
Exactly to say, in this case, Lemma~\ref{Ef_2=0} is
incorrect for $f_2$ as stated, and we could not find a substitute for $f_2$
that would provide a relevant information both for this voter model and for
the hybrid process constructed via  mixing this voter model with an exclusion 
process of any range. When the hybrid process consists of a nearest neighbor
voter model and a finite range exclusion process, an analogue of the
statement~2 i) of Theorem~\ref{pp} may 
be obtained by an appropriate, although
straightforward, modification of our arguments. Anything beyond this result 
was not possible. A reason for this was again, our failure in finding the
substitutes of  $f_1$ and $f_2$ that would  work for this
case as well as $f_1$ and $f_2$ have worked for the nearest neighbor system.

The presented above results show that, to a certain extent, the success of
our methods depend on the correct choice of the Lyapunov function.

\section*{Acknowledgements}
The authors thank FAPESP and CNPq for financial support.


\begin{thebibliography}{9}
  
\bibitem{AIM} S.~Aspandiiarov, R.~Iasnogorodsky and M.V.~Menshikov
  (1996) Passage-time moments for non-negative stochastic processes
  and an application to reflected random walks in a quadrant.  {\it
    Ann. Probab.} {\bf 24} (2), 932--960.
  
\bibitem{BCDFLS} M.~Bramson, P.~Calderoni, A.~De Masi, P.A.~Ferrari,
  J.~Lebowitz and R.~Schonmann (1986) Microscopic selection principle
  for a diffusion-reaction equation. {\it J. Statist. Phys.} {\bf
    45} (5/6), 905--920.

\bibitem{F91} C.~Cammarota and P.A.~Ferrari (1991)
\newblock Invariance principle for the branching exclusion process.
\newblock {\em Stochastic Process. Appl.} {\bf 38} (1), 1--11.

\bibitem{CD} J.T.~Cox and R.~Durrett (1995) Hybrid zones and voter
  model interfaces.  {\it Bernoulli\/} {\bf 1} (4), 343--370.
  
\bibitem{DFL} A.~De Masi, P.A.~Ferrari and J.~Lebowitz (1986)
Reaction-diffusion equations for interacting particle
systems.  {\it J. Statist. Phys.} {\bf 44} (3/4), 589--644.

\bibitem{CGLS} B.~Derrida, S.~Goldstein, J.~L.~Lebowitz and E.~Speer
  (1998) Shift equivalence of measures and the intrinsic structure of
  shocks in the asymmetric simple exclusion process. {\it J.
    Statist. Phys.} {\bf 93} (3/4), 547--571.
  
\bibitem{D88} R. Durrett (1988) {\it Lecture Notes on Particle Systems
    and Percolation.} Belmont, CA: Wadsworth.

\bibitem{D95} R. Durrett (1995) {\it Ten Lectures on Particle
    Systems.} St. Flour Lecture Notes, Lecture Notes in Math. New
    York: Springer--Verlag.


\bibitem{FMM} G.~Fayolle, V.A.~Malyshev and M.V.~Menshikov (1995) {\it
    Topics in the Constructive Theory of Countable Markov Chains.}
  Cambridge University Press.

\bibitem{FerrariP92a} P.A.~Ferrari (1992)
\newblock Shock fluctuations in asymmetric simple exclusion.
\newblock {\em Probab. Theory Related Fields\/} {\bf 91} (1), 81--101.

\bibitem{FerrariP94j} P.A.~Ferrari (1994)
\newblock Shocks in one-dimensional processes with drift.
\newblock In: G.~Grimmett (ed.), {\em Probability and Phase Transition
  (Cambridge, 1993)}, {\em NATO Adv. Sci. Ins. Ser. C, Math.
  Phys. Sci.} {\bf 420}, 35--48, Dordrecht. Kluwer Acad. Publ.

\bibitem{F} P.A.~Ferrari (1996) Growth processes on a strip.
 In: {\it
    Disordered systems (Temuco, 1991/1992)} 87--111. Travaux en Cours,
  {\bf 53}. Hermann, Paris.
  
\bibitem{FerrariP91a} P.A. Ferrari, C.~Kipnis and S.~Saada (1991)
  \newblock Microscopic structure of travelling waves in the
  asymmetric simple exclusion process.  \newblock {\em Ann Probab.}
  {\bf 19} (1), 226--244.


\bibitem{L} T.M.~Liggett (1976) Coupling the simple exclusion process.
  {\it Ann.\ Probab.} {\bf 4}, 339--356.
  
\bibitem{L85} T.M.~Liggett (1985) {\it Interacting Particle Systems.}
  Springer, Ber\-lin.

\bibitem{L99} T.M.~Liggett (1999) {\it Stochastic Interacting
  Systems. Exclusion, voter and contact processes.}
  Springer, Ber\-lin.

\bibitem{Ma} F.P.~Machado (1998) Asymptotic shape for the branching 
exclusion process. {\it Markov Processes Relat. Fields\/} {\bf 4} (4),
 535--547.
  
\bibitem{M} V.A.~Malyshev (1998) Random grammars. (Russian) 
 {\it Uspekhi Mat. Nauk\/}
{\bf 53} (2(320)), 107--134; translation in {\it Russian Math. Surveys\/}
{\bf 53}
 (2), 345--370.
  
\bibitem{MP} M.V.~Menshikov and S.Yu.~Popov (1995) Exact power estimates
  for countable Markov chains.  {\it Markov Processes Relat. Fields\/}
  {\bf 1} (1), 57--78.

\end{thebibliography}
\end{document}